\documentclass[12pt]{article}

      \usepackage{latexsym}
         \usepackage[reqno, namelimits, sumlimits]{amsmath}
         \usepackage{amssymb, amsfonts}
         \usepackage{amsthm}
 \newtheorem{theorem}{Theorem}[section]
 
 \newtheorem{definition}[theorem]{Definition}

 \newtheorem{pro}[theorem]{Proposition}

\begin{document}

\begin{center} {\Large\bf  Necessary conditions of potential  blow up for Navier-Stokes equations}\\
  \vspace{1cm}
 {\large
  G. Seregin}

 \end{center}

 \vspace{1cm}

 \vspace{1cm}
 \noindent
 {\bf Abstract }  Assuming that $T$ is a potential blow up time, we show that $H^\frac 12$-norm of the velocity field goes to $\infty$ as time $t$ approaches $T$.

 \vspace {1cm}

\noindent {\bf 1991 Mathematical subject classification (Amer.
Math. Soc.)}: 35K, 76D.


\setcounter{equation}{0}
\section{Motivation}

In the present paper, we address the following question. Consider the Cauchy problem
for the classical Navier-Stokes system
\begin{equation}\label{m1}
    \partial_t v +v \cdot \nabla v -\Delta v=-\nabla q,\qquad {\rm div} v=0
\end{equation}
in $Q_+=\mathbb R^3\times ]0,\infty[$ with the initial condition
\begin{equation}\label{m2}
    v|_{t=0}=a
\end{equation}
in $\mathbb R^3$. Here, as usual, $v$ and $q$ stand for the velocity field and for the pressure field, respectively. For simplicity, let us assume
\begin{equation}\label{m3}
    a\in C^\infty_{0,0}(\mathbb R^3)\equiv \{v\in C^\infty_0(\mathbb R^3):\,\,{\rm div} v=0\}.
\end{equation}

It is well known due to J. Leray, see \cite{Le}, that this problem has at least one weak solution obeying the  energy inequality
\begin{equation}\label{m4}
    \frac 12\int\limits_{\mathbb R^3}|v(x,t)|dx+\int\limits^t_0\int\limits_{\mathbb R^3}|\nabla v|^2dxdt'\leq \frac 12 \int\limits_{\mathbb R^3}|a|^2dx
\end{equation}
for all positive $t$. This solution is smooth for sufficiently small values of $t$. The first instant of time $T$ when singularities occur is called {\it blow up} time. By definition, $z_0=(x_0,t_0)$ is called a singular point of $v$ if it is not a regular one. The point $z_0$ is called regular if $v$ is essentially bounded in a nonempty parabolic ball  of $z_0$ \footnote{$Q(z_0,r)=B(x_0,r)\times ]t_0-r^2,t_0[$ is a parabolic ball of radius $r$ centered at the point $z_0$.}.

It is an open problem whether or not there exists
an energy solution to the Cauchy problem (\ref{m1})--(\ref{m3}) exhibiting a finite time blow up.

However, J. Leray proved some necessary conditions for $T$ to be a blow up time. They are as follows. Assume that $T$ is a  blow up time. Then, for any $3<m\leq \infty$, there exists a constant $c_m$ depending on $m$ only such that
\begin{equation}\label{m5}
    \|v(\cdot,t)\|_m\equiv \|v(\cdot,t)\|_{m,\mathbb R^3}\equiv\Big(\int\limits_{\mathbb R^3}|v(x,t)|^mdx\Big)^\frac 1m\geq \frac {c_m}{(T-t)^\frac {m-3}{2m}}
\end{equation}
for all $0<t<T$.

For the limit case $m=3$, it has been proved in \cite{ESS4} that
\begin{equation}\label{m6}
    \limsup\limits_{t\to T-0}\|v(\cdot,t)\|_3=\infty
\end{equation}
provided $T$ is a  blow up time. The interesting and open question is whether or not the following holds true:
\begin{equation}\label{m7}
  \lim\limits_{t\to T-0}\|v(\cdot,t)\|_3=\infty.
\end{equation}
The same kind of questions appears for Sobolev spaces. Coming back to the pioneer paper of J. Leray \cite{Le}, let us mention the following fact  proved there:
\begin{equation}\label{m8}
\|\nabla v(\cdot,t)\|_2\geq \frac C{(T-t)^\frac 14}
\end{equation}
for all $0<t<T$ and for some positive constant $C$ independent of $v$. This, together with the Galiardo-Nirenberg inequality and (\ref{m5}), can be extended to the following necessary condition of a finite blow up time.
\begin{equation}\label{m9}
     \lim\limits_{t\to T-0}\|v(\cdot,t)\|_{H^l}\geq \frac{\widetilde{c}_l}{(T-t)^\frac {2l-1}4}, \qquad 1/2<l<1,
     \end{equation}
     for all $0<t<T$ and for some ${\widetilde{c}_l}$, where the semi-norm  $\|\cdot\|_{H^\frac 12}$ is defined as
$$\|f\|^2_{H^l}\equiv\int\limits_{\mathbb R^3}\int\limits_{\mathbb R^3}\frac {|f(x)-f(y)|^2}{|x-y|^{3+2l}}\,dx\,dy.$$

What we are interested in here is whether or not
\begin{equation}\label{m10}
     \lim\limits_{t\to T-0}\|v(\cdot,t)\|_{H^\frac 12}=\infty.
\end{equation}
This can be regarded as the limit case in (\ref{m9}) for $l=1/2$.

Both norms $\|\cdot\|_3$ and $\|\cdot\|_{H^\frac 12}$ are very important in the mathematical theory of the Navier-Stokes equations since they are invariant with respect to their scaling:
\begin{equation}\label{m11}
    v^\lambda(x,t)=\lambda v(\lambda x,\lambda^2t),\qquad q^\lambda(x,t)=\lambda^2 q(\lambda x,\lambda^2t)
\end{equation}
for positive $\lambda$.

Actually, if the answer to the first question is positive, then (\ref{m9}) holds true by continuity of imbedding $\dot{H}^\frac 12(\mathbb R^3)$\footnote{$\dot{H}^\frac 12(\mathbb R^3)$ is the completion
of $C^\infty_0(\mathbb R^3)$ with respect to $\|\cdot\|_{H^\frac 12}$} into $L_3(\mathbb R^3)$.

In the paper, we are going to explain why (\ref{m10}) is valid. %Regarding question , an answer is still unknown.
Unfortunately, we do not know the complete answer to (\ref{m7}). Let us list some results in this direction. First of all, in \cite{S5}, a weaker version of (\ref{m7}) has been proved
$$
\lim\limits_{t\to T-0}\frac 1{T-t}\int\limits^T_t\|v(\cdot,s)\|^3_3ds=\infty.
$$
In \cite{S6}, statement  (\ref{m7}) has been demonstrated  under the additional assumption that our blow up time $T$ is of type I. More precisely, (\ref{m7}) holds true if for some $m\in ]3,\infty]$ there exists a positive constant $C_m$ depending on $m$ only such that
$$
 \|v(\cdot,t)\|_m\leq  \frac {C_m}{(T-t)^\frac {m-3}{2m}}
 $$
for all $0<t<T$. Actually, the latter condition implies the validity of it for $m=\infty$,
which means that we are dealing with blowups of type I.

Let us state our main result.
\begin{theorem}\label{mainresult} Let $v$ be an energy solution to the Cauchy problem
(\ref{m1}) and (\ref{m2}) with the smooth initial data satisfying (\ref{m3}). Let $T>0$ be a finite blow up time. Then (\ref{m10}) holds true.
\end{theorem}

Unfortunately, we still cannot justify (\ref{m7}). However, one could show the following.
\begin{pro}\label{l3} Let $v$ be an energy solution to the Cauchy problem
(\ref{m1}) and (\ref{m2}) with the smooth initial data satisfying (\ref{m3}). Let $T>0$ be a finite blow up time. Assume that for some positive number $T_1\leq T$ and for a sequence $t_k\to T_1-0$ as $k\to \infty$ the following conditions hold:
\begin{equation}\label{m12}
\sup\limits_k\|v(\cdot,t_k)\|_3=M<\infty
\end{equation}
and
\begin{equation}\label{m13}
    0<\alpha_1<\frac {t_{k+n_k}-t_k}{T-t_k}<\alpha_1<1\qquad \forall k\in \mathbb N
\end{equation}
for some nondecreasing sequence of integer $n_k$ and for real numbers $\alpha_1$ and $\alpha_2$.

Then $T_1<T$.
\end{pro}

The spirit of these two statement is that if our norms are bounded at least along a sequence
converging to a potential finite blow up time $T$, then actually $T$ is not a blow up time.
For $H^\frac 12$-norm, this is a rigorous statement while, for $L_3$-norm, it is still a plausible
conjecture. Proposition \ref{l3} says that the conjecture is true if a sequence $t_k$ converges to a potential blow up not too fast.

Now, let us shortly discuss a proof. The known way is to reduce the problem either to Liouville type theorems for bounded ancient solutions,  see \cite{KNSS} and \cite{SS3}, or to backward uniqueness for the heat operator with lower order term, see \cite{ESS4}. So far, the experience  shows  us that
working with scale-invariant norms it is preferable to utilize the second approach. Although the theory of backward uniqueness itself is relatively well understood, its realization is not an easy task and based on fine regularity results for solutions to the Navier-Stokes equations. %There are two main points to be checked. Actually, one does not need a construct an ancient solution. It is enough to construct a non-trivial solution with the help of scaling and limiting procedure which
Using the blow up technique, one can construct a non-trivial solution that
is equal to zero at the last moment of time and has a reasonable decay at infinity with respect to spatial variables. The first property easily follows from  the fact that the original solution has a finite
%$H^\frac 12$-norm and
$L_3$-norm at the blow up time $T$. And for this, boundedness along a sequence is sufficient. What is much more complicated is to construct a "blow up" solution (produced by scaling and limiting procedure) with required decay at infinity. Without such a property, the backward uniqueness  might be even wrong. The idea, how to provide such a decay, is as follows.
Applying the scaling as in \cite{S5} and \cite{S6}, we construct a local energy solution
being zero at the last moment of time $t=0$. Such kind of solutions has been introduced
 in \cite{LR1}. Here, we proceed as in \cite{KS}.
This solution has the correct decay if the initial data possess a modest decay. And this is exactly the point where the difference between $H^\frac 12$-norm and $L_3$-norm appears. We need strong compactness of initial data in $L_{2,loc}$ which is the case if one has boundedness in $H^\frac 12$-norm and is not the case if one has boundedness in $L_3$-norm. Assumption (\ref{m13}) provides required compactness at some later instance of time  for the case of $L_3$-norm.

\setcounter{equation}{0}
\section{Estimates of Scaled Solution}

There is a common part when proving  Theorem \ref{mainresult} and Proposition \ref{l3} and it is as follows. Assume that  an increasing sequence $t_k$ converges to $T$ as $k\to \infty$ and
\begin{equation}\label{sb1}
    \sup\limits_{k\in \mathbb N}\|v(\cdot,t_k)\|_{\mathcal B}=M<\infty,
\end{equation}
where $\mathcal B$ is either $H^\frac 12$ or $L_3$.  Using continuous embedding of  $\dot{H}^\frac 12$ into $L_3$ and the partial regularity theory for the Navier-Stokes equations, we can state, see similar arguments in \cite{S5} and \cite{S6}, that:
\begin{equation}\label{sb2}
    \|v(\cdot,T)\|_{3,B}<\infty.
\end{equation}
Here, $B$ is the unit ball of $\mathbb R^3$ centered at the origin.

Since $T$ is a blow up time, there exists at least one singular point at time $T$. Without loss of generality, we may assume that it is $(0,T)$.
%According to $\varepsilon$-regularity for the Navier-Stokes equations, see \cite{CKN}, there exists an universal constant such that
%\begin{equation}\label{sb2}
  %  \frac 1{r^2}\int\limits_{Q((0,T),r)}|v|^
%\end{equation}

Let us scale $v$ so that
\begin{equation}\label{sb3}
    u^{(k)}(y,s)=\lambda_kv(x,t),\quad p^{(k)}(y,s)
    =\lambda_k^2q(x,t),\quad (y,s)\in \mathbb R^3\times ]-\infty,0[,
\end{equation}
where
$$x=\lambda_ky,\qquad t=T+\lambda_k^2s,$$
$$\lambda_k=\sqrt{\frac {T-t_k}{S}}$$
and a positive parameter $S<10$ will be defined later.

Since spaces $\mathcal B$ are scale-invariant, scaled functions have the following property
\begin{equation}\label{sb4}
  \sup\limits_{k\in \mathbb N}  \|u^{(k)}(\cdot,-S)\|_{\mathcal B}=M<\infty.
\end{equation}

We know that our solution is smooth on $[-S,0[$ and we may use the uniform local estimate from \cite{LR1}
 \begin{equation}\label{sb5}
    \alpha(s)+\beta(s)\leq c\Big[\|u^{(k)}(\cdot,-S)\|_{2,{\rm unif}}+\int\limits^s_{-S}(\alpha(\tau)+\alpha^3(\tau))\,d\tau\Big],
 \end{equation}
which is valid for any $ s\in [-S,0[$ and for some positive constant $c$ independent of $k$, $s$, and $S$.  Here, the following notation have been used
$$\|f\|_{2,{\rm unif}}=\sup\limits_{x\in R^3}\|f\|_{2,B(x,r)}, \qquad \alpha(s)=\|u^{(k)}(\cdot,s)\|^2_{2,{\rm unif}},$$
$$\beta(s)=\sup\limits_{x\in R^3}\int\limits_{-S}^s\int\limits_{B(x,r)}|\nabla u^{(k)}|^2dyd\tau,$$
and $B(x,r)$ is a ball of $\mathbb R^3$ centered at the point $x$ with radius $r$.

In addition, we have the following estimate of scaled pressure, see for instance \cite{KS},
\begin{equation}\label{sb6}
\delta(0) \leq c\Big[\gamma(0)+\int\limits_{-S}^0\alpha^\frac 32(s)ds\Big]
\end{equation}
with some positive constant $c$ independent of $k$ and $S$. Here, $\gamma$ and $\delta$ are defined as
$$\gamma(s)=\sup\limits_{x\in R^3}\int\limits_{-S}^s\int\limits_{B(x,r)}| u^{(k)}|^3dyd\tau$$ and $$\delta(s)=\sup\limits_{x\in R^3}   \int\limits_{-S}^s\int\limits_{B(x,3/2)}|p^{(k)}(y,\tau)-c^{(k)}_x(\tau)|^\frac 32dy\,d\tau
$$
with some function $c^{(k)}_x\in L_\frac 32(-S,0)$ and $\gamma$ satisfies the multiplicative inequality
\begin{equation}\label{sb7}
  \gamma(s)\leq c\Big(\int\limits_{-S}^s\alpha^3(\tau)d\tau\Big)^\frac 14\Big(\beta(s)+  \int\limits_{-S}^s\alpha(\tau)d\tau\Big)^\frac 34.
\end{equation}
Now, from (\ref{sb5})-- (\ref{sb7}), it follows the existence of two positive constants $S$ and $A$ independent of $k$ such that
\begin{equation}\label{sb8}
    \sup\limits_{-S<s<0}\alpha(s)+\beta(0)+\gamma(0)+\delta(0)\leq A<\infty.
\end{equation}
And this  defines  a parameter $S$ of our scaling.

In addition to the above energy estimates, one can get bounds for higher derivatives
$\partial_tu$, $\nabla^2u$, and $\nabla p$ in $L_{\frac 98,\frac 32}$ locally. They are a simple consequence of the local regularity theory for the Stokes system. Finally, using covering by the unit balls, we can write down the following estimates
$$ \|\partial_t u\|_{L_{\frac 98,\frac 32}(B(a)\times ]-5S/6,0[)}+ \|\nabla^2 u\|_{L_{\frac 98,\frac 32}(B(a)\times ]-5S/6,0[)}$$
\begin{equation}\label{sb9}
   + \|\nabla p\|_{L_{\frac 98,\frac 32}(B(a)\times ]-5S/6,0[)}\leq C(M,a).
\end{equation}
It is worthy to note that the right hand side in (\ref{sb9}) is independent of $k$.

\setcounter{equation}{0}
\section{Limiting Procedure}

Now let us see what happens if $k\to \infty$.
Using the diagonal Cantor procedure, one can selected a subsequence, still denoted by $u^{(k)}$,
such that, for any $a>0$,
\begin{equation}\label{lp1}
   u^{(k)}\to u
\end{equation}
weakly-star in $L_\infty(-S,0;L_2(B(a)))$ and strongly in $L_3(B(a)\times ]-S,0[)$ and in
$C([\tau,0];L_\frac 98(B(a)))$ for any $-S<\tau<0$;
\begin{equation}\label{lp2}
    \nabla u^{(k)}\to \nabla u
\end{equation}
weakly in $L_2(B(a)\times ]-S,0[)$;
\begin{equation}\label{lp3}
    t\mapsto\int\limits_{B(a)}u^{(k)}(x,t)\cdot w(x) dx\to t\mapsto\int\limits_{B(a)}u(x,t)\cdot w(x) dx
\end{equation}
strongly in $C([-S,0])$ for any $w\in L_2(B(a))$;
\begin{equation}\label{lp4}
    p^{(k)}-c^{(k)}=\widetilde{p}^{(k)}\to p
\end{equation}
weakly in $L_\frac 32(B(a)\times ]-S,0[)$ for some suitable sequence $c^{(k)}\in L_\frac 32(-S,0)$.

Now, our aim is to show that $u$ is not identically zero  solution to the Navier-Stokes equations
and, moreover, it is the so-called local energy solution in the interval $]-S_1,0[$ with some $S_1\leq S$. Let us  start with the first task.

Using the inverse scaling, we have the following identity
$$\frac 1{a^2}\int\limits_{Q(a)}(|u^{(k)}|^3+|\widetilde{p}^{(k)}|^\frac 32)dy\,ds=\frac 1{(a\lambda_k)^2}\int\limits_{Q(z_T,a\lambda_k)}(|v|^3+|q-b^{(k)}|^\frac 32)dx\,dt$$
for all
$0<a<a_*=\inf\{ \sqrt{S/10},\sqrt{T/10}\}$ and for all $\lambda_k\leq 1$. Here, $z_T=(0,T)$ and $b^{(k)}(t)=\lambda_k^2c^{(k)}(s)$. Since the pair $v$ and $q-b^{(k)}$ is a suitable weak solution to the Navier-Stokes equations in $Q(z_T,\lambda_ka_*)$,  we find
\begin{equation}\label{lp5}
 \frac 1{a^2}\int\limits_{Q(a)}(|u^{(k)}|^3+|\widetilde{p}^{(k)}|^\frac 32)dy\,ds>\varepsilon
\end{equation}
for all $0<a<a_*$ with a positive universal constant $\varepsilon$.

Here, we follow arguments from \cite{S6}. Our first observation is that, by  (\ref{lp1}) and (\ref{lp4}),
\begin{equation}\label{lp6}
 \frac 1{a^2}\int\limits_{Q(a)}|u^{(k)}|^3dy\,ds\to\frac 1{a^2}\int\limits_{Q(a)}|u|^3dy\,ds
\end{equation}
for all $0<a<a_*$ and
\begin{equation}\label{lp7}
   \sup\limits_{k\in \mathbb N} \frac 1{a_*^2}\int\limits_{Q(a_*)}(|u^{(k)}|^3+|\widetilde{p}^{(k)}|^\frac 32)dy\,ds= M_1<\infty.
\end{equation}
To treat the pressure $\widetilde{p}^{(k)}$, we do the usual decomposition of it into parts. The first one is completely controlled by the pressure while the second one is a harmonic function in $B(a_*)$ for all admissible $t$. In other words,
we
have
$$\widetilde{p}^{(k)}=p^{(k)}_1+p^{(k)}_2$$
where $p^{(k)}_1$ obeys the estimate
\begin{equation}\label{lp8}
    \|p^{(k)}_1(\cdot,s)\|_{\frac 32,B(a_*)}\leq c\|u^{(k)}(\cdot,s)\|^2_{3,B(a_*)}.
\end{equation}
For the harmonic counterpart of the pressure, we have
$$\sup\limits_{y\in B(a_*/2)}  | p^{(k)}_2(y,s)|^\frac 32\leq c(a_*)\int\limits_{B(a_*)} | p^{(k)}_2(y,s)|^\frac 32dy$$
\begin{equation}\label{lp9}
\leq c(a_*) \int\limits_{B(a_*)}( |\widetilde{p}^{(k)}(y,s)|^\frac 32+|u^{(k)}(y,s)|^3)dy
\end{equation}
for all $-a_*^2<s<0$.

For any $0<a<a_*/2$,
$$\varepsilon\leq \frac 1{a^2}\int\limits_{Q(a)}( |\widetilde{p}^{(k)}|^\frac 32+|u^{(k)}|^3)dy\,ds\leq$$
$$\leq c \frac 1{a^2}\int\limits_{Q(a)}( |p_1^{(k)}|^\frac 32+|p_2^{(k)}|^\frac 32+|u^{(k)}|^3)dy\,ds\leq$$
$$\leq c \frac 1{a^2}\int\limits_{Q(a)}( |p_1^{(k)}|^\frac 32+|u^{(k)}|^3)dy\,ds+$$
$$+ca^3\frac 1{a^2}\int\limits_{-a^2}^0\sup\limits_{y\in B(a_*/2)}  | p^{(k)}_2(y,s)|^\frac 32ds.$$
From (\ref{lp7})--(\ref{lp9}), it follows that
$$\varepsilon\leq c\frac 1{a^2}\int\limits_{Q(a_*)}|u^{(k)}|^3dy\,ds+c a \int\limits_{-a^2}^0ds\int\limits_{B(a_*)}( |\widetilde{p}^{(k)}(y,s)|^\frac 32+|u^{(k)}(y,s)|^3)dy\leq$$
$$\leq c\frac 1{a^2}\int\limits_{Q(a_*)}|u^{(k)}|^3dy\,ds+ca\int\limits_{Q(a_*)}( |\widetilde{p}^{(k)}|^\frac 32+|u^{(k)}|^3)dy\,ds\leq$$
$$\leq c\frac 1{a^2}\int\limits_{Q(a_*)}|u^{(k)}|^3dy\,ds+cM_1a$$
for all $0<a<a_*/2$. After passing to the limit and picking up sufficiently small $a$, we find
\begin{equation}\label{lp10}
    0<c\varepsilon a^2\leq \int\limits_{Q(a_*)}|u|^3dy\,ds
\end{equation}
for some positive $0<a<a_*/2$. So, our limiting solution is non-trivial.

To carry on the second task, let us first recall the definition of local energy solutions

\begin{definition}\label{dlp1}
A pair of functions $u$ and $p$ defined in the space-time
cylinder $\widetilde{Q}=\mathbb R^3\times ]-S,0[$ is called a local energy weak
 Leray-Hopf solution or simply local energy solution
 to the Cauchy problem
 \begin{equation}\label{lp11}
    \partial_t u +u \cdot \nabla u -\Delta u=-\nabla p,\qquad {\rm div} u=0
\end{equation}
in $\widetilde{Q}=\mathbb R^3\times ]-S,0[$ with the initial condition
\begin{equation}\label{lp12}
    u|_{t=0}=b
\end{equation}
in $\mathbb R^3$ %Here, as usual $v$ and $q$ stand for the velocity field and for the pressure field, respectively.(\ref{d2}), (\ref{d3})
if
 the following conditions are satisfied:
 $$
u\in L_\infty(-S,0;L_{2,{\rm unif}}),\qquad\sup\limits_{x_0\in \mathbb R^3}\int\limits^0_{-S}
\int\limits_{B(x_0,1)}|\nabla u|^2dy\,ds<+\infty, $$
\begin{equation}\label{lp13} p\in L_\frac
32(-S,0;L_{\frac 32,{\rm loc}}(\mathbb R^3));\end{equation}
\begin{equation}\label{lp14}u\,\, and\,\,p\,\, meet\,\, (\ref{lp11})\,\,in
\,\,the\,\,sense\,\, of\,\,distributions;\end{equation} \begin{equation}\label{lp15} the\,
function\,\,t\mapsto\int\limits_{\mathbb R^3}u(x,t)\cdot {w}(x)\,dx
\,is\,\, continuous\,\,on\,\,[-S,0] \end{equation} for any compactly supported
function ${w}\in L_2(\mathbb R^3)$;

for any compact K, \begin{equation}\label{lp16}\|u(\cdot,t)-b(\cdot)\|_{L_{2}(K)}\to 0
\quad as\quad t\to -S+0;\end{equation}

%  to remove numbering (before each equation)
$$\nonumber \int\limits_{\mathbb R^3}\varphi|u(x,t)|^2\,dx+
2\int\limits_{-S}^{t}\int\limits_{\mathbb R^3}\varphi|\nabla
u|^2\,dxdt\leq \int\limits_{-S}^{t}\int\limits_{\mathbb
R^3}\Big(|u|^2(\partial_t\varphi+\Delta\varphi) +$$\begin{eqnarray}\label{lp17}
 +u\cdot\nabla \varphi(|u|^2+2p)\Big)\,dxdt  &&
\end{eqnarray}
for a.a. $t\in ]-S,0[$ and for all nonnegative  functions
$\varphi\in C^\infty_0(\mathbb R^3\times ]-S,1[)$;

 for any
$x_0\in\mathbb R^3$, there exists a function $c_{x_0}\in L_\frac 32
(-S,0)$ such that \begin{equation}\label{lp18}p_{x_0}(x,t)\equiv
p(x,t)-c_{x_0}(t)=p_{x_0}^1(x,t)+p_{x_0}^2(x,t),\end{equation} for $(x,t)\in
B(x_0,3/2)\times ]-S,0[$, where
$$p_{x_0}^1(x,t)=-\frac 13 |u(x,t)|^2
+\frac 1{4\pi}\int\limits_{B(x_0,2)}K(x-y): u(y,t)\otimes
u(y,t)\,dy,$$$$p_{x_0}^2(x,t)=\frac 1{4\pi}\int\limits_{\mathbb
R^3\setminus B(x_0,2)}(K(x-y)-K(x_0-y)): u(y,t)\otimes u(y,t)\,dy$$
and $K(x)=\nabla^2 (1/|x|)$.
% whenever $x\in B(x_0,3/2)$.
\end{definition}
Here, we have used the marginal Morrey space $L_{m,{\rm unif}}$ with the following finite norm
$$\|u\|_{m,{\rm unif}}=\sup\{\|u\|_{m,B(x,1)}:\,\,x\in \mathbb R^3\}.$$
Repeating arguments from \cite{KS}, we can claim that our limiting functions $u$ and $p$ satisfy
all  conditions of Definition \ref{dlp1} except condition (\ref{lp16}). In \cite{RS}, it has been shown that if $u^{(k)}(\cdot,-S)$ converges $u(\cdot,0)$ in $L_{2,{\rm loc}}$, then $u$ and $p$ is a local energy solution to Cauchy problem (\ref{lp11}), (\ref{lp12}) with $b(\cdot)=u(\cdot,-S)$.

\textsc{Proof Theorem \ref{mainresult}} By (\ref{sb4}) and by the known compactness imbedding,
we have
$$u^{(k)}(\cdot,-S) \to u(\cdot,0)$$
in $L_{2,{\rm loc}}$. Moreover, since the limit function $u(\cdot,-S)\in L_3$,
$$\|u(\cdot,-S)\|_{2, B(x,1)}\to 0$$
as $|x_0|\to\infty$. The latter, together with Theorem 1.4 from \cite{KS} and $\varepsilon$-regularity theory for the Navier-Stokes equations, gives required decay at infinity.
The last remark is that our solution has the following important property:
\begin{equation}\label{lp19}
    u(\cdot,0)=0.
\end{equation}
This follows from (\ref{sb2}) and (\ref{lp1}), see the last statement in (\ref{lp1}).
 %the fact that in any case the quantity is bounded $\|v(\cdot, t_k)\|_3$ in $k$ and by $\varepsilon$-regularity theory, $\|v(\cdot, t_k)\|_3$ is bounded as well. Then applying ,
So, we get (\ref{lp19}), for details we refer to papers \cite{S5} and \cite{S6}.
According to backward uniqueness for the Navier-Stokes, $u(\cdot,s)=0$  for any $-a_*^2<s<0$, which contradicts (\ref{lp10}). So, $z_T$ is not a singular point. Theorem \ref{mainresult} is proved.

\textsc{Proof of Proposition \ref{l3}} We still have (\ref{lp19}) but because of lack of compactness
we do not have required strong convergence in $L_{2,{\rm loc}}$.

We define $s_k$ in the following way
\begin{equation}\label{lp20}
    T+\lambda^2_k(s_k)=t_{k+n_k}=T+\lambda^2_{k+n_k}(-S).
    \end{equation}
    Then, by condition (\ref{m13}) and by (\ref{lp20}), we find that
    there exists a subsequence still denoted by $s_k$ such that
    \begin{equation}\label{lp21}
        \lim\limits_{k\to\infty}s_k =-S_0\in ]-S,0[.
    \end{equation}

Then as it follows from the last statement in (\ref{lp1})
\begin{equation}\label{lp22}
    u(\cdot,s_k) \to u(\cdot,-S_0)
\end{equation}
in $L_{2,{\rm loc}}$. Moreover, since
$$\|v(\cdot,t_k)\|_3=\|u^{(k+n_k)}(\cdot,-S)\|_3=\|u^{(k)}(\cdot,s_k)\|_3,$$
we show that
\begin{equation}\label{lp23}
    u(\cdot,-S_0)\in L_3.
    \end{equation}
So, we assume that $a_*<\sqrt{S_0/10}$ and come up with the same situation as in the proof of Theorem \ref{mainresult}, replacing $S$ with $S_0$, which means that our assumption is wrong and $z_T$ is not a singular point. Proposition \ref{l3} is proved.

\noindent\textbf{Acknowledgement} The author was partially supported by the RFFI grant
08-01-00372-a.

%\setcounter{equation}{0}
%\section{Preliminaries}

%In this section, we are going to recall the notion of local energy solutions introduced  by G.-P. %Lemarie-Rieusset, see monograph \cite{LR1} and references there. Here, we shall keep notation from  %paper \cite{KS} where a version of original statements from \cite{LR1} has been discussed.

\noindent

%G. Seregin\\
%Center for Nonlinear PDE's,\\
%Mathematical Institute, University of Oxford,UK\\
%and\\
%Steklov Institute of Mathematics at St.Petersburg, \\
%St.Peterburg, Russia

\end{document}